\documentstyle[12pt,amscd,amssymb]{amsart}

\newcommand{\Y}{\Sigma \times {\Bbb S}^1}
\renewcommand{\iff}{if and only if }
\newcommand{\Z}{\Bbb Z}
\newcommand{\R}{\Bbb R}
\newcommand{\C}{\Bbb C}
\newcommand{\mfd}{manifold with $b_1=0$}
\renewcommand{\S}{\Sigma}
\newcommand{\CP}{\Bbb C \Bbb P}
\newcommand{\SS}{{\Bbb S}^1}
\newcommand{\spinc}{$\text{Spin}^{\C}$ }
\newcommand{\A}{\cal A}
\newcommand{\B}{\cal B}
\newcommand{\inc}{\hookrightarrow}
\newcommand{\s}{\sigma}
\newcommand{\Map}{\operatorname{Map}}
\newcommand{\Diff}{\operatorname{Diff}}
\newcommand{\iso}{\cong}
\newcommand{\ar}{\rightarrow}

\theoremstyle{plain}
\begingroup
\theorembodyfont{\sl}
\newtheorem{thm}{Theorem}
\newtheorem{cor}[thm]{Corollary}
\newtheorem{lem}[thm]{Lemma}
\newtheorem{prop}[thm]{Proposition}
\endgroup

\theoremstyle{definition}
\newtheorem{defn}[thm]{Definition}

\theoremstyle{remark}
\newtheorem{rem}[thm]{Remark}

\title[Basic classes of glued manifolds]
{Constraints for Seiberg-Witten basic classes of glued manifolds}

\author{Vicente Mu\~noz}
\address{Mathematical Institute \\
24-29 St Giles'  \\
Oxford, OX1 3LB}
\email{vmunoz@@maths.ox.ac.uk}
\thanks{\hbox{$^*$}Supported by a grant from Banco de Espa\~na}
\date{June 1995}

\begin{document}

\baselineskip.55cm

\maketitle
\begin{abstract}
  We use rudiments of the Seiberg-Witten gluing theory for trivial
  circle bundles over a Riemann surface to relate de Seiberg-Witten
  basic classes of two $4$-manifolds containing Riemann surfaces of the
  same genus and self-intersection zero with those of the $4$-manifold
  resulting as a connected sum along the surface. We study examples in
  which this is enough to describe completely the basic classes.
\end{abstract}

\section{Statement of results}
\label{sec:intro}

Since their introduction nearly a year ago the Seiberg-Witten
invariants have proved to be at least as useful as their close
relatives the Donaldson invariants. These provide differentiable invariants
of a smooth $4$-manifold, whose construction is very similar in nature
to the Donaldson
invariants. Conjecturally, they give the same information about the
$4$-manifold, but they are much easier to compute in many cases, e.g.
algebraic surfaces (see~\cite{Witten}).

Problems in $4$-dimensional topology are far from solved
with these invariants. Nonetheless it is intriguing to compute them for
a general $4$-manifold. The first step towards it is obviously to
relate the invariants of a manifold with those of the manifold which
results after some particular surgery on it. Much progress has been
made~\cite{Fintushel}~\cite{Witten}. One natural case to think about
is that of connected sum along a codimension $2$ submanifold (see
Gompf~\cite{Gompf}). The typical case would be:

Let $\bar X_i$ be smooth oriented manifolds and let $\S$ be a Riemann
surface of genus $g \ge 1$.
Suppose we have embeddings
$\S \inc \bar X_i$ with image $\S_i$ representing a non-torsion
element in cohomology whose self-intersection is zero.
We form $X= \bar X_1 \#_{\S} \bar X_2$
removing tubular neighbourhoods of $\S$ in both $\bar X_i$ and gluing
the boundaries $Y$ and $\overline Y$ by some diffeomorphism $\phi$.
These boundaries are diffeomorphic to $\Y$. The diffeomorphism type of
the resulting manifold depends on the homotopy class of $\phi$. There
is an exact sequence
$$ 0 \to H^1(Y;\Z) \to H^2(X; \Z) \stackrel{\pi}{\to} G \oplus H^1(\S ;\Z)$$
with $G$ to the subgroup of $H^2(\bar X_1;\Z)/\Z[\S_1] \oplus
H^2(\bar X_2;\Z)/\Z[\S_2]$ consisting of elements $( \alpha_1, \alpha_2
)$ such that $\alpha_1 \cdot \S_1=\alpha_2 \cdot \S_2$. There are two
interpretations for this. The first one (reading the exact sequence in
homology through Poincar\'e duality) says that the $2$-homology of $X$
is composed out of the $2$-homology of $Y$ plus those cycles which
restrict to $X_1$ and $X_2$ having the same boundary $1$-cycle in $Y$
(here to be in $\pi^{-1}(G)$ means to have intersection with $Y=\Y$ a
multiple of the $\SS$ factor).
The second interpretation says that a line bundle in $X$ comes from
gluing two line bundles in $X_1$ and $X_2$
and that the possible gluings are parametrised by $H^1(Y;\Z)$.
Then the first main result we aim to prove regarding the
Seiberg-Witten basic classes is:

\begin{thm}
\label{thm:1}
Let $\bar X_i$ be smooth oriented manifolds as before and
with $b_1=0$ and $b^+>1$
and odd. Construct $X= \bar X_1 \#_{\S} \bar X_2$.
Then every basic class of $X$ lies in $\pi^{-1}(G)$. Equivalently, the
intersection of the basic class with $Y$ is $n\SS$. Moreover, $n$ is
an even integer between $-(2g-2)$ and $(2g-2)$.
\end{thm}

Our second result gives more specific information about the values of
the Seiberg-Witten invariant. It should be understood as constraints
in the possible basic classes.

\begin{thm}
\label{thm:2}
Now suppose that $\bar X_i$ are of simple type and $g \geq 2$.
Denote by $SW_X(L)$
the Seiberg-Witten invariant associated to the characteristic line
bundle $L$ (here we identify a line bundle with its first Chern class).
Fix a pair $(A_1,A_2) \in G$ such that $A_i \cdot \S_i = \pm (2g-2)$.
Then the sum $\sum_L SW_X(L)$, where $L$
runs over all characteristic line bundles in $\pi^{-1}(A_1,A_2)$
is the product of the two terms
$$
  \sum_{\{L/ L \equiv A_i \in H^2(\bar X_i)/\Z[\S_i]\} } SW_{\bar X_i}(L)
$$
Note that as $g \ge 2$ at most one of the line bundles can appear in
that sum for
each of the $\bar X_i$ because of the simple type condition.
\end{thm}

In many cases these theorems are all that we need to find out some of
the basic classes for $X$. As an example we will prove

\begin{cor}
  Suppose that either for every cycle $\gamma \in H^1(\S ; \Z)$ there exists
  a $(-1)$-embedded disc in both $X_i$ bounding $\gamma$ or
  that both $\bar X_i$ are K\"ahler manifolds and $X$ is
  deformation equivalent to $\bar X_1 \cup_{\S} \bar X_2$. Then the
  basic classes $\kappa$ of $X$ such that $\kappa \cdot \S =
  \pm (2g-2)$ are in one-to-one correspondence with
  pairs of basic classes $(\kappa_1, \kappa_2)$ for $\bar X_1$ and
  $\bar X_2$
  respectively, such that $\kappa_1 \cdot \S_1 =\kappa_2 \cdot \S_2
  =\pm (2g-2)$.
\end{cor}

{\em Acknowledgements:\/} I would like to express my gratitude to
Prof.\ Donaldson who has been an invaluable source of ideas throughout
the completion of this work. Also some discussions with
Dr.\ Kronheimer and Paul Seidel have been very illuminating. Finally I
am very grateful to Banco de Espa\~na for financial support during
these two years.

\section{Seiberg-Witten invariants}
\label{sec:2}

First of all we are recalling the basics of the Seiberg-Witten
invariants recently introduced and which are (conjecturally)
equivalent to the more classic Donaldson invariants in most relevant
cases (see~\cite{Witten}).

Let $X$ be a smooth compact oriented {\mfd } and Riemannian metric $g$,
and let a \spinc structure $c$ be given. We denote the associated
$U(2)$-bundles by $W^+=W^+_c$ and $W^-=W^-_c$ and note that Clifford
multiplication gives rise to a linear isomorphism $\rho: \Lambda^+
\otimes {\Bbb C} \to {\frak s\frak l} (W^+)$,
taking the imaginary forms into ${\frak s\frak u}(W^+)$.
Denote $L=L_c= \Lambda^2 W^+ = \Lambda^2 W^-$ the associated line bundle
of $c$. In general, the first Chern class $c_1=c_1(L_c)$ is a lift of
$w_2(X)$ to integer coefficients (what is called a characteristic
cohomology class) and for any such lift, the possible \spinc
structures with such $c_1(L)$ are parametrised by the $2$-torsion part of
$H^2(X;{\Bbb Z})$ (so if $X$ is simply connected the \spinc structures
are determined by $c_1$, which may be any characteristic class).
Said otherwise, if we fix a \spinc structure $c_0$,
the other \spinc structures are parametrised by $H^2(X;{\Bbb Z})$
with $\mu \in H^2(X;{\Bbb Z})$ corresponding to
$W_c^{\pm}=W_{c_0}^{\pm} \otimes \mu$ and $L_c = L_{c_0} \otimes \mu^2$.

The monopole equations introduced by Seiberg and Witten~\cite{Witten} for
a pair $(A, \Phi)$ of connection $A$ on the line
bundle $L$ and section $\Phi \in \Gamma (W^+)$ are
\begin{equation}
  \label{eqn:sw}
  \begin{cases}
    \not\!\!{D}_A \Phi =0 \\
    \rho(F_A^+) =(\Phi \Phi^*)_0
  \end{cases}
\end{equation}
where $\not{\!\!{D}}_A:\Gamma (W^+) \to \Gamma (W^-)$
stands for the Dirac
operator on the \spinc bundle and $(\Phi \Phi^*)_0$ is the trace-free
part of $(\Phi \Phi^*)$ interpreted as an endomorphism of $W^+$ in a
natural way.

The gauge group $\cal G = \cal C^{\infty}(X, \SS)$ acts on the
configuration space $\cal A (L_c) \times \Gamma (W^+_c)$ with quotient
$\cal B$. The moduli
space of solutions of the equations~\eqref{eqn:sw} sits in the
quotient and will be denoted by $\cal
W_{X,g}(c)$. It has the structure of a finite dimensional real
analytic space with expected dimension $d={1 \over 4}(c_1^2 - 2\chi -3
\sigma)$, where $\chi$ is the Euler characteristic of $X$ and $\sigma$
its signature.
A solution $(A, \Phi)$ is reducible (i.e. has non-trivial stabiliser)
\iff $\Phi=0$. The moduli space is always compact and, for a generic
metric $g$, is smooth away from reducibles. Whenever $b^+>0$ and
$c_1(L)$ is not torsion reducibles can
be avoided for a generic metric. Moreover,
if $b^+>1$ reducibles are not
present in generic $1$-parameter
families of metrics. Another way of obtaining a smooth moduli space is
by perturbing the equations by adding a self-dual $2$-form to
$F_A^+$~\cite{Kronheimer}.
It is also a fact that the moduli space is
orientable and that an orientation is determined by a choice of
homology orientation for $X$ (see~\cite{Witten}, \cite{Fintushel}).

\begin{defn}
  If $d<0$ then for a generic metric the moduli space is empty.
  We define the Seiberg-Witten
  invariant $SW_X(c)$ for the \spinc structure $c$ as follows: choose
  a generic metric $g$ making the moduli space smooth, compact, oriented
  (fixing a homology orientation of $X$) and of dimension the expected
  dimension $d$. Then
  \begin{enumerate}
  \item if $d<0$, put $SW_X(c)=0$.
  \item if $d=0$, let $SW_X(c)$ be the number of points of $\cal
    W_{X,g}(c)$ counted with signs.
  \item if $d>0$ and odd, we set $SW_X(c)$ to be zero. Note that this
    condition is equivalent to $b^+$ being even (since $b_1=0$).
  \item if $d>0$ and even, then we still can defined $SW_X(c)$ cutting
    down the moduli space with $d/2$ times the generator $\mu$ of the
    cohomology ring of $\cal B$ (which in this case is homotopy
    equivalent to $\CP^{\infty}$). $\mu$ also can be thought as the
    first Chern class of the $U(1)$-bundle ${\cal W}^o_{X,g}(c) \ar
    {\cal W}_{X,g}(c)$, where the superscript means framed moduli
    space.
  \end{enumerate}
\end{defn}

\begin{defn}
  One defines, for a characteristic line bundle $L$, $SW_X(L)$ to be
  the sum of $SW_X(c)$ over all \spinc structures $c$ with associated
  line bundle $L$ (note that there are a finite number of them).
\end{defn}

\begin{defn}
  Let $X$ be a compact, oriented \mfd, and suppose that $b^+ >1$. Then
  we say that a cohomology class $c_1 \in H^2(X ;\Z)$ is a {\em
  Seiberg-Witten basic class\/} (or a basic class for brevity) for $X$
  if $c_1$ is characteristic with $c_1^2 = 2\chi +3
  \sigma$ and $SW_X(c_1) \ne 0$.
\end{defn}

One important remark in place is the fact that the set of basic
classes is finite~\cite{Witten}

\begin{defn}
For $X$ compact, oriented {\mfd } with $b^+ >1$ and odd, we say that
it is of {\em (Seiberg-Witten) simple type\/} if $SW_X(L)
= 0$ whenever $d>0$.
\end{defn}

\section{Splitting along $\Y$}
\label{sec:3}

Our main interest is the study of the behaviour of the Seiberg-Witten
invariants under elementary surgeries. This amounts to splitting $X$ along
an embedded $3$-manifold $Y \subset X$. So we have $X=X_1 \cup_Y X_2$,
where $X_1$ and $X_2$ are manifolds with boundary. We orient $Y$ so
that $\partial X_1=Y$ and $\partial X_2 =\overline Y$, $Y$ with
reversed orientation.

The simplest cases are those for which $Y$ admits a metric of positive
scalar curvature. For instance, for $Y= {\Bbb S}^3$ (i.e. $X$ is a
connected sum) we have that the hypothesis $b^+(X_i)>0$ for both $X_i$
leads in a straightforward way to the vanishing of all the invariants for
$X$~\cite{Witten}. The case $b^+(X_1)=0$ and $b^+(X_2) >0$ is also of
interest and we have for instance the following theorem about the
behaviour of the Seiberg-Witten invariants under blowing-ups~\cite{Fintushel}

\begin{prop}
\label{thm:3}
If $X$ is of simple type and $\{ K_i \}$ is the set of basic classes
of $X$, then the blow-up $\tilde{X} =X \# \overline{\CP}^2$ is of simple
type and (denoting by $E$ the exceptional divisor)  the set of basic
classes is $\{ K_i \pm E \}$.
\end{prop}

The case of relevance to us is when $Y$ is a trivial circle bundle over a
Riemann surface, that is $Y=\Y$. Suppose we have an embedded $\S
\inc X$ of self-intersection $n= \S \cdot \S \geq 0$. Then we can
blow-up $X$ $n$ times in points of $\S$. Algebraically (when $X$ is
complex) this makes
perfect sense, differentiably this amounts to considering $X \#
n\overline{\Bbb C \Bbb P}^2$ and substituting the original $\S$ by its
proper transform $\S- \sum_{i=1}^n E_i$ (where $E_i$ denote the
homology class coming from the $i$-th $\overline{\Bbb C \Bbb P}^2$
summand). Therefore we can assume that $n=0$. In this case take $X_2$
to be a tubular neighbourhood of $\S$ and $X_1$ to be the closure of
the complement of it. We have a decomposition $X=X_1 \cup_Y X_2$. Here
$X_2$ is diffeomorphic to $A=\S \times D^2$.

More generally, consider the case of two manifolds $\bar{X}_1$
and $\bar{X}_2$ together with embeddings $\S \inc
\bar{X}_i$. Call the image $\S_i$ and put $n_i$ for the
self-intersection of $\S_i$. If $n_1 +n_2 \geq 0$ we can blow-up
$\bar{X}_1$ or $\bar{X}_2$ until we lower that quantity to
zero. Then we can put $N_i$ for an open tubular neighbourhood of $\S_i$,
$X_i =\bar{X}_i - N_i$, so we have an orientation reversing
diffeomorphism $\phi$ between the boundaries of $X_1$ and $X_2$. This
lets us
construct $X=X_1 \cup_{\phi} X_2$ which we call {\em connected sum of
  $\bar{X}_1$
and $\bar{X}_2$ along $\S$\/} (obviously the diffeomorphism type of
$X$ depends on the homotopy class of $\phi$). A simple extension of the
arguments in~\cite{Gompf} gives the following

\begin{prop}
  The diffeomorphism type does not depend on which points we blow-up.
  More concretely, if we blow-up $X_i$ in $s_i$ points with $s_1 +s_2
  = n_1+n_2$ and we choose a diffeomorphism $\phi$, there is a unique
  homotopy class of diffeomorphisms $\phi'$ such that blowing-up $X_i$
  in $s'_i$ points with $s'_1 +s'_2
  = n_1+n_2$ the result is diffeomorphic to the former.
\end{prop}

Therefore the only choice (up to diffeomorphism) involved in all
this process is the identification between the boundaries (that is an
element of $\pi_0(\Diff^-(Y))$). Note that when $n_1 \geq 0$ and $n_2
\geq 0$ we can lower both quantities to zero to have $Y=\Y$.

Suppose now that $b_1=0$ for both $\bar X_1$ and $\bar X_2$.
Consider embeddings $\S \inc X_i$ with images $\S_i$. If $[\S_i]
\in H^2(\bar X_i; \Z)$
is a non-torsion element, then the cohomology exact sequence
for the pair $(\bar X_i,X_i)$ gives
$$  H^1(X_i;\Z) \iso H^1(\bar X_i;\Z) $$
$$  H^3(X_i;\Z) \iso H^3(\bar X_i;\Z) $$
$$  0 \to H^2(\bar X_i;\Z)/\Z[\S_i] \to H^2(X_i;\Z) \to
H_1(\S ;\Z) \to 0$$
where the last map is the composition
$$H^2(X_i;\Z) \to H_2(X_i,
\partial X_i;\Z) \stackrel{\partial}{\to} H_1(\Y ;\Z)
\stackrel{i_*}{\to} H_1(\S \times D^2 ;Z)=H_1(\S;\Z)$$
Hence there is a (non-canonical) splitting $H^2(X_i;\Z) =
H^2(\bar X_i;\Z)/\Z[\S_i] \oplus H_1(\S ;\Z)$.

\noindent The Mayer-Vietoris sequence for $X=X_1 \cup X_2$ gives (we
drop $\Z$ in the notation)
$$H^1(X )= H^1(X_1) \oplus H^1(X_2)$$
$$0 \to H^1(Y) \to H^2(X) \to H^2(X_1) \oplus H^2(X_2) \to
H^2(Y) \iso H^1(\S) \oplus H^2(\S)$$
So $b_1(X)=0$. Under the splitting above, we can describe the last map as
\begin{eqnarray*}
  (H^2(\bar X_1)/\Z[\S_1] \oplus H_1(\S)) \oplus (H^2(\bar
  X_2)/\Z[\S_2] \oplus H_1(\S )) & \to & H^1(\S;\Z) \oplus
  H^2(\S;\Z) \\
  ( \alpha_1, \beta_1, \alpha_2, \beta_2) & \mapsto & (\beta_1 -\beta_2,
  (\alpha_1 -\alpha_2) \cdot \S)
\end{eqnarray*}
(the identifications between homology and cohomology groups are
through Poincar\'e duality).
Calling $G$ to the subgroup of $H^2(\bar X_1;\Z)/\Z[\S_1] \oplus
H^2(\bar X_2;\Z)/\Z[\S_2]$ consisting of elements $( \alpha_1, \alpha_2
)$ such that $\alpha_1 \cdot \S_1=\alpha_2 \cdot \S_2$ (note that
these pairings make sense), we have that
\begin{equation}
  \label{eqn:coh}
  0 \to H^1(Y;\Z) \to H^2(X; \Z) \stackrel{\pi}{\to} G \oplus H^1(\S ;\Z)
\end{equation}

\begin{rem}
  Equation~\eqref{eqn:coh} is the most explicit description we can get
  of the picture. This admits two different interpretations. If we
  think in terms of the homology, the first term corresponds to the
  $2$-homology of
  $Y$, i.e. $\S$ and those classes of the form $\gamma \times \SS$,
  where $\gamma \in H_1(\S;\Z)$. $\pi^{-1}(G)$ are the classes obtained
  by gluing cycles coming in $\bar X_1$ with cycles in $\bar X_2$
  intersecting $\S$ on the same number of points. This process is not
  well-defined and actually choosing different representatives of given
  homology classes in both $\bar X_i$ we can get all of $\pi^{-1}(G)$.
  The preimage of $ H^1(\S ;\Z)$ corresponds to $2$-cycles which
  intersection with $Y$ is a $1$-cycle contained in $\S$ (that is,
  cycles with part in $X_1$ and in $X_2$ going through the neck). This last bit
  is not canonically defined as explained above. \linebreak
  If we think in terms of line bundles and their first Chern classes,
  $\pi(L)$ is the restriction of $L$ to the two open manifolds $X_i$
  and $H^1(Y;\Z)$ express the different ways in which two line bundles
  on each of $X_i$ could be glued to give a line bundle on $X$.
\end{rem}

\begin{rem}
  The characteristic numbers are related as follows:
  $$ \chi_X= \chi_{\bar X_1} + \chi_{\bar X_2} +4g-4$$
  $$ \sigma_X = \sigma_{\bar X_1}+\sigma_{\bar X_2}$$
  Therefore $b^+$ and $b^-$ are both increased by $2g-1$ and $2\chi
  +3\sigma$ by $8g-8$.
\end{rem}

\begin{rem}
\label{rem:1}
  Call $m_i$ the minimal number of $(\S_i \cdot D_i)$, for $D_i \in H^2(X_i
  ;\Z)$. Put $m$ for the least common multiple of $m_1$ and $m_2$ and
  $d=m_1m_2/m$ for the greatest common divisor. Then $m$ is the
  minimal number of $(\S \cdot D)$, for $D \in H^2(X ;\Z)$ and
  the map $\pi$ above is surjective if and only if $m_1$ and $m_2$ are
  coprime. In general the cokernel is isomorphic to $\Z/d\Z$.
\end{rem}

\section{Seiberg-Witten equations for $\Y$}
\label{sec:4}

During the last years there has been a great deal of work on
developing the gluing theory for computing the Donaldson invariants of
a manifold $X=X_1 \cup_Y X_2$ out of information from $X_1$ and $X_2$
(see~\cite{Austin}, \cite{Braam}).
The standard technique is to pull apart $X_1$ and $X_2$ so that we are led
to consider metrics giving $X_i$ a cylindrical end and $L^2$-solutions
of the equations in these open manifolds.
This process has an analogue in the Seiberg-Witten setting, first
introduced in~\cite{Kronheimer}. We refer there for the notations used
here. The analogue of the Chern-Simons
functional is ${1 \over 4\pi^2}C$ on $\B_3=\A_3 \times \Gamma(W_3)/
\operatorname{Map} (Y, \SS)$ taking values on $\R/\Z$ and given by
$$C(A, \Phi)= \int (A-B)\wedge F_A + \int <\Phi, \not{\!\partial}_A \Phi>$$
where $B \in \A_3$ is a fixed connection on $Y$. The downward gradient
flow equations for this functional are
\begin{equation}
  \begin{cases}
    dA/dt = -* F_A + i \tau(\Phi, \Phi) \\ d\Phi/dt =-
    \not{\!\partial}_A \Phi
  \end{cases}
\end{equation}
with $\tau: \bar W_3 \times W_3 \to \Lambda^1 (Y)$ from
Clifford multiplication and $\not{\!\!\partial}_A$ the $3$-dimensional
Dirac operator coupled with $A$. The critical points correspond to translation
invariant solutions on the tube. So they are solutions to
\begin{equation}
  \label{eqn:sw3d}
  \begin{cases}
    * F_A = i \tau(\Phi, \Phi) \\  \not{\!\partial}_A \Phi=0
  \end{cases}
\end{equation}
The reducible solutions are those for which $\Phi=0$ and therefore
$F_A=0$ and $c_1(L)=0$.

Now we pass on to study the equations~\eqref{eqn:sw3d} for $Y=\Y$. We
choose for $Y$ a metric which is rotation invariant.
Let $L \to \Y$ be a line bundle. We can pull it back to $\S \times
[0,1]$ under the map identifying $\S \times \{0\}$ with $\S \times
\{1\}$ and denote it by $L$ again. Then topologically this line bundle is
the pull back of some line bundle $M$ on $\S$. Clearly $c_1(M)$ is
the restriction of $c_1(L)$ to $\S$ and $L$ is obtained from the pull
back of $M$ to
$\S \times [0,1]$ by gluing with some gauge transformation $g \in
\Map(\S, \SS)$. The homotopy class of $g$ is the component of $c_1(L)$
in $H^1(\S;\Z) \otimes H^1(\SS;\Z) $
under the isomorphisms $[\S ; \SS ] \iso H^1(\S;\Z) \iso H^1(\S;\Z)
\otimes H^1(\SS;\Z)$, where the last one is product with the
fundamental class of the circle. Now denote by $\theta$ the coordinate
in the $\SS$-directions. Any connection $A$ on $L$ has a representative in
its gauge equivalence class with no $d\theta$ component. This is
unique up to constant gauge
transformation (i.e. a gauge pulled-back from $\S$). So giving $A$ (up
to gauge) is
equivalent to giving a family $A_{\theta}$, $\theta \in [0,1]$ of
connections on $M$ (up to
constant gauge) with the condition $A_1 = g^*(A_0)$, with $g$ in the
homotopy class determined by $L$.

The \spinc structure on $Y$ restricts to a \spinc structure
on $\S$ and therefore $W_3$ restricted to $\S$ decomposes as $( \Lambda^0
\oplus \Lambda^{0,1})
\otimes \mu$ for some line bundle $\mu$ on $\S$. Note that $M= K_{\S}
\otimes \mu^2$, where $K_{\S}$ is the canonical bundle of $\S$.
The Dirac operator $\not\!\partial_A : \Gamma(W_3) \to \Gamma(W_3)$ is
identified with $i{\partial \over \partial \theta} + \sqrt{2}(\bar
\partial_{A_{\theta}} +\bar \partial_{A_{\theta}}^*)$.
Recall that the connection
$A_{\theta}$ is on the bundle $M$ and induces uniquely a connection on
the bundle $\mu$.

\begin{lem}
  Let $(A, \Phi) \in \A_3 \times \Gamma(W_3)$. Consider the family
  $A_{\theta}$, $\theta \in [0,1]$ determined by $A$ and write $\Phi =
  (\s_1, \s_2) \in \Gamma (W_3) = \Omega^0 (\mu) \oplus  \Omega^{0,1}
  (\mu)$. Then the solutions of~\eqref{eqn:sw3d} correspond to
  solutions of:
  \begin{equation} \label{eqn:sw2d}
    \begin{cases}
      {\partial \s_1 \over \partial \theta }= -  \sqrt{2}i \, \bar
        \partial_{A_{\theta}}^* \s_2 \\
      {\partial \s_2 \over \partial \theta }= - \sqrt{2}i \, \bar
        \partial_{A_{\theta}} \s_1 \\
      i{\partial A_{\theta} \over \partial \theta }= \s_1 \s_2^* +\s_2
        \s_1^* \\
      -2 \, i \Lambda F_{A_{\theta}} = | \s_1 |^2 - |\s_2|^2
    \end{cases}
  \end{equation}
\end{lem}

In the third equation $\s_1 \s_2^* +\s_2 \s_1^* \in \Omega^1$ is a
real form. Recall that
the connection $A_{\theta}$ is equivalent to the holomorphic structure
$\bar \partial_{A_{\theta}}$, so we can rewrite that line as
$${\partial \over \partial \theta } (\bar \partial_{A_{\theta}})=
-i\,\s_1^* \s_2$$

%\begin{rem}
%  The system of equations~\eqref{eqn:sw2d} is overdetermined. Writing
%  $F_{A_{\theta}}= \partial_{A_{\theta}} \bar \partial_{A_{\theta}} + \bar
%  \partial_{A_{\theta}} \partial_{A_{\theta}}$ and using that
%  $${\partial \over \partial \theta } (\partial_{A_{\theta}})= - i\,\s_1
%  \s_2^*$$
%  we can take derivatives with respect $\theta$
%  to get $-2i\,\Lambda F_{A_{\theta}} - | \s_1 |^2 + |\s_2|^2$
%  constant. Therefore the fourth equation in~\eqref{eqn:sw2d} holds
%  for every $\theta$ whenever it does for $\theta =0$.
%\end{rem}

\begin{prop}
\label{prop:1}
  Let $(\s_1, \s_2) \in \Omega^0 (\mu) \oplus  \Omega^{0,1}
  (\mu)$ and $A_{\theta}$, $\theta \in [0,1]$ such that
  \begin{equation} \label{eqn:sw2d'}
    \begin{cases}
      {\partial \s_1 \over \partial \theta }=  - \sqrt{2}i \,\bar
        \partial_{A_{\theta}}^* \s_2 \\
      {\partial \s_2 \over \partial \theta }= - \sqrt{2}i \,\bar
        \partial_{A_{\theta}} \s_1 \\
      {\partial \over \partial \theta } (\bar \partial_{A_{\theta}})=-i\,\s_1^*
      \s_2
    \end{cases}
  \end{equation}
  Then $A_{\theta}$, $\s_1$ and $\s_2$ are constant and either $\s_1 =
  0$ and $\bar \partial_{A_0}^* \s_2=0 $ or $\s_2 =0$ and $\bar
  \partial_{A_{0}} \s_1=0$.
\end{prop}

\begin{pf}
$$
  \bar \partial ^* \bar \partial \s_1 = - {1 \over  \sqrt{2}i} \bar
  \partial ^*( {\partial
  \s_2 \over \partial \theta}) = {\partial \over \partial \theta}
  (-{1 \over  \sqrt{2}i}\bar \partial ^* \s_2)+ {1 \over  \sqrt{2}i}
  {\partial \over \partial \theta}
  (\bar \partial ^*) \s_2= {\partial ^2 \s_1 \over \partial \theta^2} -
  {1 \over  \sqrt{2}} \s_1 |\s_2|^2
$$
where we drop subindices by convenience of notation.
First integrate along $\S$ by parts to get for every $\theta$
$$ \int_{\S} |\bar \partial \s_1|^2 - \int_{\S} <{\partial ^2 \s_1
  \over \partial \theta^2} , \s_1 >
  +{1 \over  \sqrt{2}} \int_{\S} | \s_1^*\s_2|^2 =0 $$
This equation makes sense in $\SS$ and we can integrate again by parts
to get
$$||\bar \partial \s_1||^2 +||{\partial \over \partial
  \theta} \s_1||^2 +{1 \over  \sqrt{2}} || \s_1^*\s_2||^2 =0 $$
The result is immediate from this.
\end{pf}

\begin{cor}
\label{cor:0}
  If the line bundle $L$ admits any solution to~\eqref{eqn:sw3d} then
  $L$ is pulled-back from $\S$ and any solution is invariant under
  rotations in the $\SS$ factor.
\end{cor}

\begin{rem}
  We can paraphrase corollary~\ref{cor:0} saying that any basic class
  is orthogonal to $H_1(\S;\Z) \otimes H_1(\SS ;\Z) \inc H_2(Y;\Z)
  \inc H_2(X;\Z)$. This is in fact a very natural phenomenon, for if
  we put $T_{\gamma} = \gamma \times \SS$ then $T_{\gamma}$ is a torus
  of self-intersection zero and hence $K \cdot T_{\gamma} = 0$ for any
  basic class $K$.
\end{rem}

Now let $L$ be a characteristic line bundle which is the pull-back of
a line bundle in $\S$. Since $\S \cdot \S =0$ we have that $c_1(L)
\cdot \S$ is
even. Consider a \spinc structure with associated line bundle $L$
(this is the pull-back of a \spinc structure on $\S$ of the same
kind). The description of solutions of the Seiberg-Witten equations in
algebraic varieties~\cite{Kronheimer2} \cite{Okonek} gives an
identification of the moduli space of solutions of
equations~\eqref{eqn:sw3d} with the space of vortices.

\begin{cor}
\label{cor:1}
  Let $2k=(2g-2) - |c_1(L) \cdot \S|$. Then the moduli space of
  solutions of the translation invariant equations~\eqref{eqn:sw3d} on
  $\Y$ is empty if $k<0$, $s^k (\S)$ whenever $c_1(L) \cdot \S \ne 0$
  (here $s^k(\S)$
  stands for the symmetric power of $\S$). If $c_1(L) \cdot \S =0$ then
  the moduli space consists uniquely of reducibles and is isomorphic
  to the Jacobian of line bundles of degree $g-1$ over $\S$.
\end{cor}

\begin{rem}
  When $m=c_1(L) \cdot \S =0$ we have to perturb the
  equations~\eqref{eqn:sw} by adding a self-dual form. Choosing a form
  invariant under translations and rotations, this produces
  the effect of perturbing the last equation in~\eqref{eqn:sw2d} by
  adding a small $2$-form to the curvature. The proof of
  proposition~\ref{prop:1} goes through and the methods for proving
  corollary~\ref{cor:1} give a moduli space of solutions $s^{g-1}(\S)$.
\end{rem}

Theorem~\ref{thm:1} is a consequence of corollary~\ref{cor:1} and the following

\begin{cor}
  Let $c$ be a \spinc structure for $X$ with associated line bundle
  $L$. If $c$ is a basic class for $X$ then $L|_Y$ is pulled-back from
  $\S$.
\end{cor}

\section{Gluing theory}
\label{sec:5}

In the
Seiberg-Witten context there is a parallel of the usual Floer theory
for the Donaldson invariants~\cite{Braam}. Some nice few remarks about
the case relevant to us appear in~\cite{Donaldson}.
Generally for a three-manifold $Y$ and a line bundle $L|_Y$ on $Y$ (we
use this as a matter of notation as $L|_Y$ will be the restriction of a
line bundle $L$ in a $4$-manifold containing $Y$), we perturb the
equations for the translation invariant solutions of~\eqref{eqn:sw3d} on
the tube until the solutions are finite and non-degenerate. Then one
defines
Seiberg-Witten Floer homology groups $HFSW_*(Y;L|_Y)$ as the Floer groups
in the usual context. The cohomology groups
$HFSW^*(Y;L|_Y)$ are naturally identified with the homology groups of $Y$ with
reversed orientation, hence giving a natural pairing between them.

In the case of $Y=\Y$, we need to consider characteristic line bundles
$L$ whose
restrictions to $Y$ have $c_1(L|_Y)= 2m \SS$, for $|m| \leq g-1$, as
already established in theorem~\ref{thm:1}.
Fix $m$, i.e. the topological type of $L|_Y$ and put $k=(g-1) -m$.
The picture is analogue to the instanton
theory. Instead of using a perturbation of the
equations~\eqref{eqn:sw3d},
we can work with the moduli space of solutions itself,
which by corollary~\ref{cor:1} is $M_{\S} = s^k(\S)$, whenever $m \neq 0$
(when $M_{\S}$ consists only of irreducibles). When $m=0$ we perturb
the equations as explained before to be in the same situation. So
$HFSW_*(Y;L)$ is identified with the
homology of $M_{\S}$.

Let $X_1$ be an open manifold with cylindrical end $Y$. For a \spinc
structure over $X_1$ whose associated line
bundle $L$ is of the required type, the limit
values of solutions to the Seiberg-Witten equations give the element
$\Phi_{X_1}(c)$ of $HFSW_*(Y;L)$. As before $\Phi_{X_1}(L)$
will denote the sum over all possible $c$ giving rise to the same
$L$.
When we have two open manifolds $X_1$ and $X_2$ which we want to glue
along the common boundary $Y$ (with a fixed diffeomorphism of the
boundaries), we have in general an indeterminacy for
choosing the identification of the line bundles over $Y$ resulting in
different line bundles for $X= X_1 \cup_Y X_2$.
In the case of section~\ref{sec:3}, i.e. when $b_1(X_i)=b_1(\bar X_i)=0$, the
possibilities are parametrised by $H^1(Y;\Z)$.
In this case we have to refine the groups $HFSW_*(Y;L|_Y)$ to keep track
of the homotopy class of this identification.
For that, we mod out the space of solutions by the component of the
identity of $\Map(Y,\SS)$, instead of using all of it (i.e. we work
with $\tilde{\B}_3=\A_3 \times \Gamma(W_3)/
\operatorname{Map} (Y, \SS)_o$). The result is a
number of copies of $M_{\S}$ parametrised by $H^1(Y;\Z)$, which is
called $\tilde M_{\S} = \sum_{i \in H^1(Y)} M_{\S}^{(i)}$,
producing groups  $HFSW_*\tilde{\ }(Y;L|_Y) \iso \sum_{i \in H^1(Y)}
H_*(M_{\S}^{(i)})$.
Since $b_1=0$, the invariants lift to $\tilde{\Phi}_{X_1}(L)$.
Clearly there is an addition map
$$ HFSW_*\tilde{\ }(Y;L|_Y) \stackrel{a}{\to} HFSW_*(Y;L|_Y) $$
which recuperates the original invariant, i.e.
$a(\tilde{\Phi}_{X_1}(L))=\Phi_{X_1}(L)$.

Now when we glue solutions coming from $X_1$ to solutions from $X_2$,
the first thing to have in mind is that the copies of $M_{\S}$ in
which both of them live determine uniquely a gluing of the line
bundles over the boundary (and hence the line bundle $L$ on $X$).
For instance, if $(A, \phi)$ is a solution of the Seiberg-Witten
equations in $X$ for the line bundle $L$, which splits as $(A_1,\Phi_1)$ with
limit $(a,\phi) \in  M_{\S}^{(i)}$
and  $(A_2,\Phi_2)$ with limit $(a,\phi) \in  M_{\S}^{(j)}$, then $i-j
\in H^1(Y)$ determines $L$ and $i+t$ and $j+t$ will determine the same
$L$ for any
$t \in H^1(Y ;\Z)$ (this can be thought as
transferring $t$ from $X_2$ to $X_1$ through the neck).
When we pull $X_1$ and $X_2$ apart introducing metric with larger and
larger tube lengths, every solution on $X$ to the
equations~\eqref{eqn:sw} appears as solutions in $X_1$ and $X_2$ giving
the same boundary value and such that the determined line bundle is
the one we started with. So the pairing between $\Phi_{X_1}(L|_{X_1})$
and $\Phi_{X_2}(L|_{X_2})$ corresponds to
using all possible gluings of $L|_Y$.

\begin{prop}
\label{thm:a}
  Let $L$ be a line bundle over $X=X_1 \cup_Y X_2$. Then the pairing
  $$\tilde{\Phi}_{X_1}(L|_{X_1}) \cdot \tilde{\Phi}_{X_2}(L|_{X_2})$$
  in $HFSW_*\tilde{\ }(Y;L|_Y)$ is equal to the summation
  (see~\eqref{eqn:coh} for definition of $\pi$)
  $$ \sum_{\{L' / \pi(c_1(L')) =\pi(c_1(L)) \} } SW_X(L')$$
  i.e. over all $L'$ whose restrictions to both $X_i$ are isomorphic
  to the ones of $L$.
\end{prop}

\begin{rem}
In the case of $X_2=A= D^2 \times \S$ we cannot do the same thing as
we have to take into account the homology of $A$, which is $H^1(\S
;\Z)$.
Therefore we only can lift the Floer groups to $\tilde M_{\S} / H^1(\S
;\Z)$, which is formed by copies of $M_{\S}$ parametrised by $\Z
[\S]$. We call this $M'_{\S}= \sum_{n \in \Z[\S]}  M_{\S}^{(n)}$.
\end{rem}

We also have a forgetful map $$HFSW_*\tilde{\ }(Y;L) \stackrel{\alpha}{\to}
HFSW'_*(Y;L) = HFSW_*\tilde{\ }(Y;L)/ H^1(\S;\Z)$$

\begin{prop}
\label{thm:b}
  Let $L$ be a line bundle over $\bar X_1 = X_1 \cup_Y A$. Then the
  pairing
  $$\alpha(\tilde{\Phi}_{X_1}(L|_{X_1}))  \cdot \Phi'_A(L|_A)$$
  in $HFSW'_*(Y;L|_Y)$ is equal to the summation
  $$ \sum_{\{L' / c_1(L') \equiv c_1(L) \pmod{\S} \} } SW_{\bar X_1}(L')$$
  Furthermore, if $\bar X_1$ is of simple type and $m \neq 0$, then
  at most one of the $L'$ can appear in the sum above, since at most
  one of them has $c_1(L)^2 = 2 \chi +3\sigma$.
\end{prop}

Theorem~\ref{thm:2} is a consequence of the following

\begin{thm}
\label{thm:c}
  Now suppose that $\bar X_i$ are of simple type and $g \geq 2$.
  Fix a pair $(\kappa_1,\kappa_2) \in G$ such that $\kappa_i \cdot
  \S_i = \pm (2g-2)$ and $\kappa_i$ are characteristic.
  Then the sum $\sum_{L \in \pi^{-1}(\kappa_1,\kappa_2)} SW_X(L)$
  is the product of the two terms
  $$
  \sum_{\{L/ c_1(L) \equiv \kappa_i \in H^2(\bar X_i)/\Z[\S_i]\} }
  SW_{\bar X_i}(L)=SW_{\bar X_i}(\kappa_i)
  $$
\end{thm}

\begin{pf}
  In the case $c_1(L) = \pm (2g-2) \SS$, one has $k=0$ so $M_{\S}$
  is a point and $H_*(M_{\S}) \iso \Z$. Now the solutions for $A=\S
  \times D^2$ are all pulled-back from $\S$. So $\Phi_A(L|_A)$
  consists of a $1$ in one $\Z$ and zero in the rest and hence
  $a(\Phi_A(L|_A)) =1$. Now
  proposition~\ref{thm:b} says
  $$ n_{\kappa_i} =SW_{\bar X_i}(\kappa_i)
  =\alpha(\tilde{\Phi}_{X_i}(L|_{X_i}))  \cdot \Phi'_A(L|_A) =
  {\Phi}_{X_i}(L|_{X_i})$$
  Now the results comes from applying proposition~\ref{thm:a}.
\end{pf}

\begin{rem}
  Note that the same argument is not applicable to the case $k>0$ as
  plugging in $A$ we only get (knowing $SW_{\bar X_i}(L)$) the values
  of ${\Phi}_{X_i}(L|_{X_i})$
  in $H_0(M_{\S}^{(j)})$, but we need the higher homology. One expects
  that no new basic classes are going to appear from pairs
  $(\kappa_1,\kappa_2) \in G$ such that $|\kappa_i \cdot
  \S_i| < 2g-2$.
\end{rem}

\section{Examples}
\label{sec:6}

This section is devoted to examples in which the information already
gathered in propositions~\ref{thm:a} and~\ref{thm:b} is enough to find the
basic classes of the glued manifold.

\begin{thm}
\label{thm:6.1}
  Suppose that for a basis of homology cycles $\gamma \in H^1(\S ;
  \Z)$ there are embedded $(-1)$-discs in both $X_i$ bounding some
  embedded $1$-cycle representing $\gamma$, then the
  basic classes $\kappa$ of $X$ such that $\kappa \cdot \S =
  \pm (2g-2)$ are in one-to-one correspondence with
  pairs of basic classes $(\kappa_1, \kappa_2)$ for $\bar X_1$ and
  $\bar X_2$
  respectively, such that $\kappa_1 \cdot \S_1 =\kappa_2 \cdot \S_2
  =\pm (2g-2)$. Moreover, $\kappa$ is determined in an explicit way.
\end{thm}

\begin{pf}
  By a $(-1)$-embedded disc as above we understand an embedding
  $(D^2, \partial D^2) \hookrightarrow (X_i, \partial X_i)$ with the
  boundary going to an embedded curve representing $\gamma$ and
  such that a {\em small\/}
  perturbation has only one point of intersection with the original
  disc, lying in the interior and with sign $-1$.
  So when the situation of the theorem is given, we can glue together
  the discs (in a generally not unique way) to get an embedded sphere
  of self intersection $-2$, say $D_{\gamma}$.
  Now call $T_{\gamma}$ to the torus $\gamma \times \SS \inc \Y$. We
  have the obvious fact $T_{\beta} \cdot T_{\gamma}=0$ and $D_{\beta}
  \cdot T_{\gamma} = (\beta \cdot \gamma)$. So there is a torus of
  self-intersection zero and a sphere of self-intersection $-2$
  intersecting in one point. This implies~\cite{Fintushel} that $X$ is
  of simple type and that the basic classes vanish on all of these
  cohomology classes.
  What is more, the $T_{\gamma}$ and $D_{\gamma}$ generate a primitive
  sublattice
  $V \subset H^2(X ;\Z)$ hence $V \oplus V^{\perp} =H^2(X ;\Z)$ and
  \begin{equation} \label{eqn:seq5}
    0 \ar \Z[\S] \ar V^{\perp} \stackrel{\pi}{\ar} G
  \end{equation}
  Now let $\kappa$ be a basic class for $X$ such that $\kappa \cdot \S
  = 2g-2$. We have argued that
  $\kappa \in V^{\perp}$, theorem~\ref{thm:1} tells us how the image
  of $\kappa$ under $\pi$ is. So there are basic classes $\kappa_1$
  and $\kappa_2$ in $X_1$ and $X_2$ such that $\kappa \cdot \S
  =\kappa_1 \cdot \S_1  = \kappa_2 \cdot \S_2 =2g-2$.
  For $g \geq 2$, we have that $\kappa^2 \neq (\kappa +n \S)^2$ for $n
  \neq 0$, so at most one of the $\kappa +n \S$ can be basic class.
  Thus the statement of the theorem.
\end{pf}

\begin{rem}
  To come to terms with theorem~\ref{thm:6.1} we need to split up the
  sequence~\eqref{eqn:seq5}. Suppose $\S_i$ are both primitive. In the
  words of remark~\ref{rem:1}, $m_1=m_2=m=1$. We choose cohomology
  classes $D \in H^2(X ;\Z)$, $D_1 \in H^2(\bar X_1 ;\Z)$, $D_2 \in
  H^2(\bar X_2
  ;\Z)$ with $D \cdot \S =D_1 \cdot \S_1 = D_2 \cdot \S_2=1$.
  Let $W$ be the primitive sublattice in $V^{\perp}$
  containing $\S$ and $D$, and
  analogously for $W_1$ and $W_2$. Then $\pi: W^{\perp} \iso
  W_1^{\perp} \oplus W_2^{\perp}$ and $W/ \Z [\S] \iso \Delta
  \subset W_1/\Z[\S_1] \oplus W_2/\Z[\S_2]$, taking $[D]$ to $([D_1],[D_2])$.
  Also $D$ provides an
  isomorphism $W \iso W/\Z [\S] \oplus \Z [\S]$ and the same for
  $D_1$ and $D_2$, so $$H^2(\bar X_i) \iso W_i^{\perp} \oplus
  W_i/\Z[\S_i] \oplus \Z[\S_i]$$ $$H^2(X) \iso V \oplus W_1^{\perp}
  \oplus W_2^{\perp} \oplus \Delta \oplus \Z [\S]$$
  Then if we choose $D^2 =D_1^2 +D_2^2$ (easy to arrange), we have (in
  the decomposition above)
  that for $\kappa_i = \alpha_i +(2g-2) D_i +r_i \S_i$ basic classes
  for $X_i$, the corresponding basic class for $X$ is $\kappa = 0 +
  \alpha_1 + \alpha_2 + (2g-2)D + (r_1 +r_2 +2)\S$ (the coefficient of
  $\S$ is found out using the fact that $\kappa^2 = 2\chi +3 \sigma$). Formally
  \begin{equation}
    \kappa = \kappa_1 + \kappa_2 + 2 \S
  \end{equation}
\end{rem}

\begin{rem}
  When both $\bar X_i$ are symplectic manifolds and $\S_i$ are
  symplectic submanifolds, Gompf~\cite{Gompf}
  has proved that $X$ can be given a symplectic structure (regardless
  of the homotopy class of the chosen gluing $\phi$). Taubes~\cite{Taubes}
  \cite{Taubes2}
  has proved that the canonical class $K=-c_1(TX)$ is a basic class
  and that for any other basic class $\kappa \neq \pm K$, one has
  $|\kappa \dot [\omega]| <  K \dot [\omega]$. Since, in the
  notation of the last theorem, $T_{\gamma} \dot [\omega] =0$, none
  of the $K + \sum n_{\gamma}T_{\gamma}$ can be basic classes unless
  all $n_{\gamma}=0$. Then in the formula of proposition~\ref{thm:a} only
  one term appears in either side. Notice that Taubes also proves that
  this number is $\pm 1$.
\end{rem}

\begin{rem}
  The result of the last remark falls very short since it does not
  even tell us about the other basic classes that might appear when we
  glue two basic classes $K_i$ for $\bar X_i$ with $K_i \cdot \S
  =2g-2$ but $K_i$ are not the canonical classes.
  In some situations we get more: suppose that (both) $\bar X_i$ have come
  out as the blow-up of some symplectic manifolds $M_i$ at points on
  $\S'_i \inc M_i$ (and $\S_i$ is the proper transform of $\S'_i$)
  and that the cohomology class defined by the
  symplectic forms in both $M_i$ where Poincar\'e dual to $[\S'_i]$
  (i.e. $[\S'_i]$ are ample classes).
  Then one has $|\kappa_i \dot [\S_i]| <  K \dot [\S_i] =2g-2$ for
  every basic class $\kappa_i$ in $\bar X_i$. So we conclude that the
  only basic classes with $\kappa \cdot \S = \pm (2g-2)$ for $X$ are
  $\kappa = \pm K$.
\end{rem}

When the manifolds involved are complex surfaces and $\S_i$ are
embedded complex curves, there is a
preferred identification between the boundaries of $X_i$ given by the
only $\phi$ which identifies holomorphically the holomorphic normal
bundles of $\S_i$ in $\bar X_i$.

\begin{prop}
  Let $\cal Z \stackrel{\pi}{\to} \Delta$ be a family of complex
  surfaces. Suppose that $Z_t= \pi^{-1}(t)$ are smooth for $t \ne
  0$ and that $X=Z_0$ is the union of two surfaces $\bar{X}_1$
  and $\bar{X}_2$ intersecting in a normal crossing
  (see~\cite{Friedman}). Then the diffeomorphism type of $X$ is
  obtained by the necessary blow-ups and a connected sum of
  $\bar{X}_1$ and $\bar{X}_2$ along $\S$ with the preferred
  identification alluded above.
\end{prop}

\begin{pf}
  To prove this result first we construct another deformation family which
  general member is $Z_t$ for $t \ne 0$ but the fibre over $0$ is the
  union of the blow-ups of $\bar X_i$ at $s_i$ arbitrary points on
  $\S_i$ (with
  $s_1 +s_2 = n_1 +n_2$). Without loss of generality, we can suppose
  the case of one point in $\bar X_1$. The problem is local, so we
  pick a small (Zariski) neighbourhood $\cal U$ of the point in $\cal
  Z$ such that $U_t= Z_t \cap \cal U$ is embedded
  in $\C^3$. The intersection $C_t = U_t \cap U_0$ is a reducible
  curve which can be written as $C_{t,1} \cup C_{t,2}$ with $C_{t,i}
  \subset \bar X_i$ (reducing even further $\cal U$ we can suppose
  that the only intersection of $ C_{t,i}$ with $\S$ is the given
  point and that this intersection is transverse). Consider
  $\cup_t \ U_t \times \{t\}$ and its divisor $\cup_t \ C_{t,2} \times \{t\}$.
  We blow-up such divisor, which does not affect any of $U_t$, $t \ne
  0$ nor $\bar X_2$ and has the result of blowing-up $\bar X_1$ at the
  given point. Now glue this new (Zariski) open affine to the
  complement of the point in $\cal Z$. We have the required
  deformation \linebreak
  For completing the proof we note that when $n_1+n_2 =0$ the
  deformation of $\bar X_1 \cup \bar X_2$ is diffeomorphic to the
  connected sum along $\S$ alluded in the statement.
\end{pf}

\begin{thm}
  Suppose that both $\bar X_i$ are K\"ahler manifolds and that $X$ is
  deformation equivalent to $\bar X_1 \cup_{\S} \bar X_2$. Then the
  basic classes $\kappa$ of $X$ such that $\kappa \cdot \S = \pm
  (2g-2)$ are in one-to-one correspondence with
  pairs of basic classes $(\kappa_1, \kappa_2)$ for $X_1$ and $X_2$
  respectively, such that $\kappa_1 \cdot \S_1 =\kappa_2 \cdot \S_2
  =\pm (2g-2)$.
\end{thm}

\begin{pf}
  Firstly, it is known after Witten~\cite{Witten} that all K\"ahler
  manifolds with $b^+>1$ are of simple type.
  As in the proof of theorem~\ref{thm:6.1} we just need to prove that if
  $\kappa$ is a basic class for $X$ and $T= \sum n_{\beta}T_{\beta} \neq
  0$ in $H^1(Y; \Z)$ then $\kappa +T$ is not basic class. In the
  K\"ahler case we know that the basic classes are in $H^{1,1}$, so it
  is enough to show that $T \notin H^{1,1}$. But $T^2=0$ and $T \dot
  [\omega] =0$, for the symplectic form $\omega$. If $T$ were in
  $H^{1,1} \cap H^2(X;\Z)$, it would represent a divisor with $T^2=0$
  and orthogonal to an ample class, but this is impossible.
\end{pf}

\begin{rem}
  Last two theorems combined with proposition~\ref{thm:3} give information
  on the basic classes of a K\"ahler manifold which can be deformed to
  an algebraic variety with normal crossings, knowing the basic
  classes of the irreducible components.
\end{rem}

\end{document}